\documentclass[11pt]{article}
\usepackage{epic,latexsym,amssymb}
\usepackage{amsfonts}
\usepackage{amscd}
\usepackage{amsmath}
\usepackage{graphicx}
\usepackage{color}
\usepackage{caption,
caption}
\usepackage{tikz}

\usepackage{float}

\newtheorem{thm}{Theorem}

\newtheorem{prop}[thm]{Proposition}
\newtheorem{lem}[thm]{Lemma}

\newtheorem{problem}{Problem}
\newtheorem{defn}{Definition}

\newcommand{\sub}{\mathrm{sub}}

\newcommand{\proof}{\noindent\textbf{Proof. }}

\newcommand{\qed}{$\Box$}

\newcommand{\QEDmark}{\mbox{\textsc{qed}}}
\newcommand{\proofStarter}[1]{\textsc{#1} }

\def\vertex(#1){\put(#1){\circle*{2}}}
\def\vertexo(#1){\put(#1){\circle{2}}}
\def\vert(#1){\put(#1){\circle*{1.5}}}
\def\verto(#1){\put(#1){\circle{1.5}}}
\def\lab(#1)#2{\put(#1){\makebox(0,0)[c]{#2}}}
\definecolor{DarkGreen}{rgb}{0.2, 0.6, 0.3}

\definecolor{electricindigo}{rgb}{0.44, 0.0, 1.0}

\begin{document}

\title{A note on sub-total domination in graphs}
\author{$^{1,2}$Randy Davila\\
\\
$^1$Department of Pure and Applied Mathematics\\
University of Johannesburg \\
Auckland Park 2006, South Africa \\
\\
$^2$Department of Mathematics \\
Texas State University \\
San Marcos, TX 78666, USA \\
\small {\tt Email: rrd32@txstate.edu}
}

\date{}
\maketitle

\begin{abstract}
Let $G$ be a simple and finite graph without isolated vertices. In this note we study a degree sequence derived invariant called the \emph{sub-total domination number}, denoted $\sub_t(G)$. This invariant originally appeared in \cite{total Slater} and serves as a lower bound on $\gamma_t(G)$, where $\gamma_t(G)$ denotes the heavily studied \emph{total domination number} of $G$.  
\end{abstract}


\section{Introduction}
Domination in graphs is widely studied and a heavily applied notion in graph theory. Indeed, domination and its variants and generalizations appear in vast quantities in the mathematical literature; see for example \cite{CaroRoditty, Delavina,  Favaron, DomBook1, DomBook2, Adriana, MHAYbookTD, Pepper, Rautenbach}. Of the many variants of domination, total domination is arguably one of the most natural. Given a graph $G$, and a set of vertices $S$ in $G$, $S$ is a \emph{total dominating set} if every vertex in $G$ has a neighbor in $S$. The minimum cardinality of a total dominating set in $G$ is the \emph{total domination number} of $G$, denoted by $\gamma_t(G)$. It is well known that determining the total domination number of a general graph is in the class of $NP$-complete decision problems \cite{NP}, and as such, a significant amount of research has been devoted to finding easily computable upper and lower bounds on $\gamma_t(G)$; see for example the monograph \cite{MHAYbookTD} which details and surveys total domination. 

As previously mentioned, finding computationally efficient bounds on $\gamma_t(G)$ is desired. However, in a much more general fashion, it is of great interest to find computationally efficient bounds for any $NP$-hard graph invariant. With this in mind, we make note that the degree sequence of a graph has been shown to yield such desired bounds. Two well known examples are the \emph{residue} and the \emph{annihilation number} of a graph, which serve as respective lower and upper bounds on the computationally difficult \emph{independence number} of a graph \cite{Favaron residue,PepperThesis}. With regards to domination, the lesser known degree sequence derived invariants known as the \emph{slater number} and the \emph{sub-k-domination number} serve as respective lower bounds on the \emph{domination number} and \emph{k-domination number} of a graph \cite{sub_k_dom, Slater}. We remark that these degree sequence results are special cases of the recently introduced \emph{degree sequence index strategy} (DSI-strategy) \cite{DSI}.  
 
\medskip
\noindent\textbf{Definitions and Notation.}
All graphs in this paper will be considered finite simple graphs without isolated vertices. Let $G = (V,E)$ be a graph. We will denote the order and size of $G$ by $n = n(G) = |V(G)|$ and $m = m(G) = |E(G)|$, respectively. When the dependence on $G$ is clear, we will write $n$ in place of $n(G)$. Two vertices $v,w\in V(G)$ are said to be neighbors if $vw\in E(G)$. The open neighborhood of $v\in V(G)$, denoted by $N_G(v)$, is the set of neighbors of $v$, whereas the close neighborhood of $v$ is the set $N_G[v] = N_G(v)\cup \{v\}$. The degree of $v\in V(G)$ is the cardinality of $N_G(v)$, and will be denoted by $d_G(v)$. The maximum and minimum vertex degree among all vertices of $G$ will be denoted by $\Delta(G)$ and $\delta(G)$, respectively. A graph $G$ is called $k$-regular if $d_G(v) = k$ for all $v\in V(G)$. A regular graph is a graph that is $k$-regular for some integer $k\ge 0$. 

The degree sequence of $G$, is the sequence consisting of the vertex degrees in $G$ listed in non-increasing order, and will be denoted $D(G) = \{\Delta(G) = d_1, \dots , d_n = \delta(G)\}$. For brevity, we may write the number of vertices realizing each degree in superscript. For example, the path $P_n$, on $n$ vertices, may have degree sequence written $D(P_n) = \{2^{n-2}, 1^{2}\}$. If a sequence of non-negative integers $D$ has the property that $D = D(G)$, for some graph $G$, then we say that $D$ is a graphic sequence, and that $D$ is realizable by $G$. We note that a given graphic sequence may have more than one graph which realizes $D$. 

A set of vertices $S\subseteq V(G)$ is a total dominating if every vertex in $G$ has a neighbor in $S$, and such a set will be called a TD-set of $G$. The cardinality of a smallest TD-set in $G$ is the total domination number of $G$, denoted by $\gamma_t(G)$, and such a set will be called a $\gamma_t(G)$-set. For other graph terminology and definitions, we will follow \cite{MHAYbookTD}. 

We will also make use of the notation $[k] = \{1, \dots, k\}$.

\section{Sub-total domination}
In this section we present our main results. First we recall the definition of the sub-total domination number, originally defined in \cite{total Slater}, and denoted $sl_t(G)$. Keeping our notation and terminology consistent with \cite{sub_k_dom}, we will use $\sub_t(G)$ in place of $sl_t(G)$.  
\begin{defn}
If $G$ is an isolate-free graph with order $n$ and degree sequence $D(G) = \{\Delta(G) = d_1,\dots, d_n = \delta(G)\}$, the sub-total domination number $\sub_t(G)$, is defined as the smallest integer $k$ such that $\sum_{i=1}^{k}d_i \ge n$.
\end{defn}
 
With the definition of sub-total domination now defined, we remark that $\sub_t(G)$ can be computed in $O(n)$ time. Because of the simplicity of computing $\sub_t(G)$, and the difficulty of computing $\gamma_t(G)$, the following theorem serves as one of our main results. We remark that this theorem first appeared in \cite{total Slater} without proof. 
\begin{thm}[\cite{total Slater}] \label{sub-total domination}
If $G$ is an isolate-free graph, then
\begin{equation*}
\gamma_t(G)\ge \sub_t(G),
\end{equation*}
and this bound is sharp.
\end{thm}

\proof
Let $G$ be a graph with order $n$, degree sequence $D(G) = \{\Delta(G) = d_1, \dots, d_n = \delta(G)\}$, and $S$ be a $\gamma_t(G)$-set. Next, we order the vertices of $S$, $s_1, \dots, s_{|S|}$, so that $d_G(s_1) \ge \dots \ge d_G(s_{|S|})$. By definition, every vertex is totally dominated by a vertex in $S$; that is, every vertex has a neighbor in $S$. Thus, $V(G) = \cup_{v\in S}N_G(v)$, which implies, 
\begin{equation*}
n = \Big |\bigcup_{v\in S}N_G(v)\Big | \le \sum_{v\in S}|N_G(v)| = \sum_{i=1}^{|S|}d_G(s_i). 
\end{equation*}
In particular, we have established, 
\begin{equation*}
\sum_{i=1}^{|S|}d_G(s_i) \ge n.
\end{equation*} 
Next observe that the $i$-th term of $D(G)$ is greater than or equal to the $i$-the degree of the list of vertices from $S$, and thus, we have the following inequality,
\begin{equation*}
\sum_{i=1}^{|S|}d_{i} \ge \sum_{i=1}^{|S|}d_G(s_i) \ge n.
\end{equation*}
That is,
\begin{equation}\label{one}
\sum_{i=1}^{|S|}d_{i} \ge n.
\end{equation}
Since $\sub_t(G)$ is the smallest integer satisfying (\ref{one}), it follows that $\gamma_t(G) = |S|\ge \sub_t(G)$, and the lower bound has been proven. 

To see that this bound is sharp, consider the star $K_{1,n-1}$ on $n\ge 2$ vertices. Then, $\gamma_t(K_{1,n-1}) = 2$, and $\sub_t(K_{1,n-1}) = 2$. 
\qed 

Theorem \ref{sub-total domination} is sharp for non-trivial stars. However, stars are a special case of a more general concept. Namely, if $G$ is a connected graph with order $n \ge 2$ and maximum degree $\Delta(G) = n - 1$, then choosing a maximum degree vertex and an arbitrary neighbor of this vertex forms a TD-set, and hence, $\gamma_t(G) = 2$. Moreover, the highest vertex degree summed with the next highest vertex degree will be greater than $n$, and so $\sub_t(G) = 2$. In particular, since no vertex of $G$ will have degree $n$, it follows that $\sub_t(G) \ge 2$. We combine these ideas with the following proposition.  
\begin{prop}\label{trivialOb1}
If $G$ is a connected graph with order $n\ge 2$ and maximum degree $\Delta(G) = n - 1$, then $\gamma_t(G) = \sub_t(G) = 2$. 
\end{prop}

There exists graphs $G$ for which $\gamma_t(G) = 2$ and $\Delta(G) \neq \Delta(G) - 1$. Double stars (trees with exactly two non-leaf vertices) are one such example. With this in mind, we next generalize Proposition \ref{trivialOb1} to a statement on graphs $G$ with $\gamma_t(G) = 2$. That is, since $\sub_t(G) \ge 2$, we obtain the following proposition. 
\begin{prop}\label{trivialOb2}
If $G$ is an isolate-free graph with $\gamma_t(G) = 2$, then $\gamma_t(G) = \sub_t(G)$. 
\end{prop}

A simple lower bound on the total domination number of isolate-free graphs can be found by dividing the order by the maximum degree, see Chapter 2, Theorem 2.11. in \cite{MHAYbookTD}. With the following theorem we show that the sub-total domination number improves on this bound.  
\begin{thm}\label{trivial}
If $G$ is an isolate-free graph with order $n$ and maximum degree $\Delta(G)$, then $\gamma_t(G) \ge \sub_t(G) \ge n/ \Delta(G)$.
\end{thm}
\proof 
Let $G$ be an isolate-free graph with order $n$ and maximum degree $\Delta(G)$. The left hand side of the inequality is a restatement of Theorem \ref{sub-total domination}.  Thus, in order to prove this result, it suffices to show $\sub_t(G) \ge n / \Delta(G)$. By definition, we have
\begin{equation*}
\sum_{i=1}^{\sub_t(G)}d_i \ge n.
\end{equation*}
Next observe that $\Delta(G)\ge d_i$ for each $i \in [\sub_t(G)]$, and thus
\begin{equation*}
\sub_t(G)  \Delta(G) = \sum_{i=1}^{\sub_t(G)}\Delta(G) \ge\sum_{i=1}^{\sub_t(G)}d_i \ge n,
\end{equation*}
Hence, $\sub_t(G) \ge n/ \Delta(G)$, and the proof of the theorem is complete. \qed

\section{Properties of $\sub_t(G)$} 
In this section we provide various fundamental properties of the sub-total domination number. We begin with a closed formula for $\sub_t(G)$ in the case that $G$ isolate-free and $k$-regular.  
\begin{prop}\label{regular graphs}
If $k\ge 1$ is an integer and $G$ is a $k$-regular graph with order $n$, then $\sub_t(G) = \lceil n / k \rceil$. 
\end{prop}
\proof
Let $k\ge 1$ be an integer, and let $G$ be a $k$-regular isolate-free graph with order $n$. By definition of sub-total domination, we have
\begin{equation*}
\sub_t(G) k= \sum_{i=1}^{\sub_t(G)} k \ge n.
\end{equation*}
It follows that $\sub_t(G) \ge n / k$. Since $\sub_t(G)$ is the smallest integer satisfying this inequality, we obtain $\sub_t(G) = \lceil n/ k \rceil$, and the proof of the proposition is complete. \qed

Next we consider sub-total domination of disjoint isolate-free graphs. In particular, we show that sub-total domination is subadditive with respect to disjoint unions of graphs. 
\begin{lem}\label{Disjoint Union}
If $G$ and $H$ are isolate-free graphs, then $\sub_t(G) + \sub_t(H) \ge \sub_t(G \cup H)$.
\end{lem}
\proof
Let $G$ and $H$ be disjoint graphs with degree sequences $D(G) = \{ \Delta(G) = d_{1}^{G}, \dots, d_{n_1}^{G} = \delta(G)\}$ and $D(H) = \{ \Delta(H) = d_{1}^{H}, \dots, d_{n_2}^{H} = \delta(H)\}$. By definition of sub-total domination, we have
\begin{equation*}
\sum_{i=1}^{\sub_t(G)}d_{i}^{G}\ge n_1,
\end{equation*}
and, 
\begin{equation*}
\sum_{i=1}^{\sub_t(H)}d_{i}^{H}\ge n_2. 
\end{equation*}
Thus, 
\begin{equation*}
\sum_{i=1}^{\sub_t(G)}d_{i}^{G} + \sum_{i=1}^{\sub_t(H)}d_{i}^{H} \ge n_1+n_2 = n(G\cup H). 
\end{equation*}
Denote the degree sequence of $G\cup H$ by $D(G\cup H) = \{\Delta(G\cup H) = d_{1}^{*}, \dots, d_{n_1+n_2}^{*} = \delta(G\cup H)\}$. Since degree sequences are listed in non-increasing order, it follows that
 \begin{equation*}
\sum_{i=1}^{\sub_t(G)+\sub_t(H)}d_{i}^{*} \ge \sum_{i=1}^{\sub_t(G)}d_{i}^{G} + \sum_{i=1}^{\sub_t(H)}d_{i}^{H} \ge n_1 + n_2 = n(G\cup H). 
\end{equation*}
That is, 
 \begin{equation}\label{three}
\sum_{i=1}^{\sub_t(G) + \sub_t(H)}d_{i}^{*} \ge n(G\cup H). 
\end{equation}
Since $\sub_t(G\cup H)$ is the smallest integer satisfying (\ref{three}), it follows that $\sub_t(G\cup H) \le \sub_t(G) + \sub_t(H)$, and the proof of the lemma is complete. \qed

It is easy to see that the total domination number is additive with respect to unions of disjoint graphs; that is, for disjoint isolate-free graphs $G$ and $H$, $\gamma_t(G\cup H) = \gamma_t(G) + \gamma_t(H)$. With this in mind, the following theorem serves as an improvement on Theorem \ref{sub-total domination} when considering the union of disjoint graphs.  
\begin{thm}
If $G$ and $H$ are isolate-free graphs, then 
\begin{equation*}
\gamma_t(G\cup H) \ge \sub_t(G) + \sub_t(H) \ge \sub_t(G\cup H).
\end{equation*} 
\end{thm}
\proof Let $G$ and $H$ be isolate-free graphs. By Lemma \ref{Disjoint Union} $\sub_t(G) + \sub_t(H) \ge \sub_t(G\cup H)$. Moreover, since total domination is additive with respect to disjoint unions, $\gamma_t(G\cup H) = \gamma_t(G) + \gamma_t(H)$. By Theorem \ref{sub-total domination}, $\gamma_t(G) \ge \sub_t(G)$ and $\gamma_t(H) \ge \sub_t(H)$. Thus, $\gamma_t(G\cup H) \ge \sub_t(G) + \sub_t(H)$, and the theorem is proven. \qed


\section{Conclusion and Open Problems}
In this note we have studied fundamental properties of $\sub_t(G)$. However, we have not studied many classes of graphs for which $\gamma_t(G) = \sub_t(G)$. Since $\sub_t(G)$ is easily computable, we suggest the following problem.
\begin{problem}\label{prob1}
Characterize all graphs $G$ for which $\gamma_t(G) = \sub_t(G)$. 
\end{problem}

Problem \ref{prob1} is surely difficult, and leads to the question of asking if determining a graph $G$ satisfies $\gamma_t(G) = \sub_t(G)$ is $NP$-complete. The analogous question for sub-domination and domination is is known to be $NP$-complete \cite{total Slater}, and so this provides evidence that this may indeed be the case. 

There exists many lower bounds on the total domination number of a graph, and it remains to be shown how sub-total domination compares with most of these bounds. Thus, we further suggest the following problem.
\begin{problem}
Compare $\sub_t(G)$ with known lower bounds on $\gamma_t(G)$. 
\end{problem}

\end{document}